\documentclass[11pt]{amsart}
\usepackage{epsfig}
\usepackage{latexsym}
\usepackage[all]{xy}
\usepackage{amsmath}
\usepackage{amsfonts}
\usepackage{amssymb}
\usepackage{xypic}

\def\R{\mathbb R}

\def\Z{\mathbb Z}

\def\S{\mathbb S}

\def\I{\mathbb I}

\def\t{\triangle}
\def\bn{\bigskip\noindent}
\newtheorem{Def}{Definition}[section]
\newtheorem{Thm}{Theorem}[section]

\newtheorem{Lemma}{Lemma}
\newtheorem{Prop}{Proposition}[section]

\usepackage{epsfig}
\def\scfig #1 #2 {\resizebox{#2}{!}{\includegraphics{#1}}}
\usepackage{fullpage}

\usepackage{xypic}
\usepackage{psfrag}
\begin{document}
\title{\sc Topological Colored Tverberg Theorem\\ and the Reduction Lemma}
\author{Satya Deo}
\address{Harish-chandra Research Institute (HRI)  Allahabad}
\email{ sdeo@hri.res.in }
\date{}

\maketitle

\let\thefootnote\relax\footnotetext{Author supported by the Senior Scientist Platinum Jubilee Fellowship of the National Academy of Sciences, India (2013)} .\\

\begin{abstract} In this paper we present a proof of the BMZ Reduction Lemma with a motivational perspective, and state this lemma for maps to manifolds using the classical definition of cohomological dimension. The lemma proved  and utilized in \cite{bmz} gives a geometrical insight in the proof of the BMZ colored topological Tverberg theorem.
\end{abstract}

\section{Introduction}
Let us recall the well-known topological Tverberg theorem  ~\cite{bar-shl-sz}:
{\Thm  Let $d\geq 1, \;r\geq 2$ a prime and  $N=(d+1)(r-1)$. Consider the $N$-simplex $\triangle_N$. Suppose $f: \t_N \to \R^d$ is a continuous map. Then there exists a family $\sigma_1,\sigma_2,\cdots\sigma_r$ of $r$ disjoint faces of $\t_N$ such that} $$\bigcap^{r}_{1}f(\sigma_i)\neq\phi.$$

\bn The above Theorem is a generalization of the classical Tverberg Theorem~\cite{tve} which says that if we have $N=(d+1)(r-1)+1$ points in $\R^d$ in general position, then we can partition them into $r$ disjoint subsets $F_1, F_2,\cdots , F_r$ such that the convex hulls of all these $F_i$'s have a non-empty intersection. The classical Tverberg Theorem has been proved for all positive integers $r\geq 1$ whereas the topological Tverberg Theorem is proved only when $r$ is a prime power. When $r$ is not a prime power, a counterexample to the topological Tverberg theorem has recently been published (see \cite {fri}).

\bn What would be the colored version of the topological Tverberg Theorem ?
This question was first studied by B$\acute{a}$r$\acute{a}$ny, F$\ddot{u}$redi and Lov\'{a}sz~\cite{bar-fur-lov}, who proved that if we take 21 points in the plane in general position such that 7 are red, 7 are blue and 7 are  green, then we can decompose them into 7 triangles, each triangle having vertices of different colors, such that their intersection is non empty. Later, B$\acute{a}$r$\acute{a}$ny and Larman~\cite{bar} studied the above question further and proved the following general result : If we take $3r$ points in the plane in general position where $r$ are of red color, $r$ are of blue color and  $r$ are of green color, then we can find $r$ disjoint  triangles, each triangle having vertices of different colors, such that their intersection will be non empty. For higher dimensional Euclidean spaces $\R^d$, B$\acute{a}$r$\acute{a}$ny and Larman asked the question:  Given $r$ and $N\geq (d+1)r$ points in $\R^d$, determine the smallest $t$ such that  if we take $d+1$ color classes of size $\leq t$, we can find  a family of $r$ disjoint $d$-simplices having vertices of different colors, such that their intersection is non empty.  $\check{Z}$ivaljevi$\acute{c}$ and Vre$\acute{c}$ica~\cite{vre-ziv} introduced the concept of chessboard complexex and showed that when $r$ is a prime, the colored Tverberg Theorem holds for $t\geq 2r-1$, and also for any $r\geq 2$ due to Bertrand's postulate that there is a prime between $r$ and $2r$.  This colored Tverberg Theorem of ZV attracted a lot of attention, but could not be considered very satisfactory for two reasons. The first reason was that the classical Tverberg Theorem does not appear as a special case of this Theorem, and the second reason was that several colored points were left out unaccounted. In 2009, Blagojevi$\acute{c}$, Matschke and Ziegler(see ~\cite{bmz}~\cite{bmz2}~\cite{bmz3}) formulated and proved a remarkable new colored topological Tverberg Theorem which is stated as follows:

{\Thm (BMZ)  Let $d\geq 1,\; r\geq 2$ prime and $N=(d+1)(r-1)$. Consider the $N$-simplex $\t_N$ whose vertices are colored into disjoint color classes
$$C_0,C_1,\cdots ,C_m\;\;\;m\geq d+2$$ such that for each $i$, the size $|C_i|$ of these colored classes is at most $r-1$. Then given any continuous map $f: \t_N\to\R^d$ there exists a family of $r$ disjoint rainbow faces $\sigma_1,\sigma_2, \cdots ,\sigma_r$ of $\t_N$ such that
$$f(\sigma_1)\cap\cdots\cap f(\sigma_r)\neq \phi.$$
}

Here  a face $\sigma_i$ is said to be a {\bf rainbow} face if each vertex of $\sigma_i$ is of different color, i.e., $|\sigma_i\cap C_j|\leq 1\;\;\forall\; i,\; 1\leq i\leq r$ and $\forall\; j,\;\;1\leq j\leq m$.

\bn The above Theorem, to be called the BMZ Theorem in the sequel, turns out to be very interesting topological colored Tverberg Theorem since the classical Tverberg Theorem becomes a special case  of this Theorem when the size of each color class is one, and also since all points of $\t_N$ are accounted for.
In this  paper we present  a simple minded proof of the BMZ Reduction Lemma (see~\cite{bmz} page 3, Reduction of Theorem 2.1 to Theorem 2.2) with a motivational perspective. We also state the BMZ reduction Lemma for maps to manifolds using the classical definition of cohomological dimension over $\Z$ rather than the one employed in (\cite{bmz3}, Lemma 2.1). We remark here that BMZ Theorem for manifolds has been proved when $r$ is prime.  The corresponding Theorem when $r$ is a prime power is not resolved yet, though the topological Tverberg Theorem has already been proved for any prime power $r$. However, the Reduction Lemma is true for any $r$, prime or otherwise.

\section{Reduction Lemma}
Let $d\geq 1,\;r\geq 2$ and $N=(d+1)(r-1)$. Color the vertices of the simplex $\t_N$ by $m$ colors, $m\geq d+2$. Suppose the color classes are denoted by $C_0,C_1\cdots ,C_m$.  If $|C_i|\leq r-1,\;\forall \; i=1,2,\cdots ,m$ then we say that
$C_0,C_1\cdots ,C_m$ is a {\bf general coloring} of $\t_N$.  Since $|C_i|\leq r-1$ for each $i$, there has to be at least $d+2$ colors.
 Among all the general colorings of $\t_N$, let us consider the colorings
$$C_0,C_1\cdots ,C_{d+1}$$ where $|C_i|= r-1$ for every $i$ except one color class whose size is 1. This accounts for all the points of $\t_N$. Without loss of generality we may occasionally consider the last color class viz.,
$C_{d+1}$, to be of size 1. Such colorings of $\t_N$ will be called the {\bf special colorings} of $\t_N.$

\bn {\bf Note :} It may be noted that any general coloring $G_0,G_1\cdots ,G_m,\;\;m\geq d+1$, of the vertices of $\t_N$ can be obtained from a given special coloring $$C_0,C_1\cdots ,C_{d+1}$$ by removing certain number of vertices from each color class and creating new color classes (one or more) in finitely many steps.
It is assumed that when we remove the single vertex of $C_{d+1}$, this singleton class disappears creating  another class of the same color with possibly bigger size.

\bn In this process  while  the number of vertices of $\t_N,\;N=(d+1)(r-1)$ remains the same, the class sizes become smaller and smaller and the number of color classes becomes larger and larger. In the limiting case when the size of each  color class becomes 1,  the number of color classes will become $N+1=(d+1)(r-1)+1$.

\bn In terms of the above definitions we can state the Blogojevi$\acute{c}$-Matschke-Ziegler colored Theorem as follows (see \cite{bmz} Theorem 2.1):

{\Thm {\bf ( BMZ Theorem)}  Let $d\geq 1,\; r\geq 2$ prime and $N=(d+1)(r-1)$  and $C_0,C_1\cdots ,C_m$ be a general coloring  of $\t_N$. Then given any continuous  map $f:\t_N\to\R^d$ there exists a family  of $r$ disjoint rainbow  faces $\sigma_1, \sigma_2,\cdots , \sigma_r$ of $\t_N$ such that
$$ f(\sigma_i)\cap\cdots\cap f(\sigma_r)\neq\phi.$$}

{\Lemma {\bf (Reduction Lemma)}  Let $d\geq 1,\; r\geq 2$  and $N=(d+1)(r-1).$  If the BMZ Theorem is true for special colorings for the parameters $(d+k, r, \R^{d+k}),\;k\geq 0,$ then the BMZ theorem is true for all general colorings of $\t_N$ for the parameters $(d,r,\R^d)$. }

\bn  We now give a detailed proof of the above Lemma following a simple minded approach with a motivational perspective.

\bn Let us start with a special coloring  $C_0,\cdots ,C_d, C_{d+1}$ of $\t_N,\;\;N=(d+1)(r-1),$ and a continuous map  $f:\t_N\to\R^d$ and assume BMZ for special colorings. First we prove the following :

{\Prop  Let us remove one vertex from some color class, say $C_0$, and give it a new color creating an additional color class, say $C_{d+2}$ such that $|C_{d+2}|=1$. Then BMZ theorem is true for the parameters $(d,r,\R^d)$ for the coloring  $C_0, C_1\cdots , C_{d+1},C_{d+2};\;\; |C_0|=r-2, \;|C_i|=r-1, \; |C_{d+1}|=1=|C_{d+2}|$.}

\begin{proof}
To prove this, we add one new vertex to $C_0$ of the same color as that of  $C_0$, making it $C_0'$, add additional $r-2$ vertices of the color of $C_{d+2}$ to make it $C_{d+2}'$ so that $|C_{d+2}'|=r-1$ (see Fig 1). Thus we have the coloring
$$ C_0', C_1',\cdots , C_{d+1}', C_{d+2}',\;\;C_i\subset C_i'\;\;\;\forall \; i$$
and $|C_i'|=r-1\;\;\forall\;i\neq d+1,\;|C_{d+1}'|=1$. Now we consider a new simplex $\t_{N'},\;N'=(d+2)(r-1)$ and regard $\t_N$ as the front face of $\t_{N'}$. Thus, there are $r-1$ additional vertices in $\Delta_{N'}-\Delta_N$.

$$\bordermatrix{ & & &  &  & \cr
& & &  & &  \cr
 C_{d+1} & \bullet &   &\cr
  C_{d} & \bullet &  \bullet &  \bullet  &  \bullet &  \bullet  \cr
                & \vdots & \vdots & \ddots & \vdots & \vdots\cr
                C_{1} & \bullet &  \bullet &  \bullet  &  \bullet &  \bullet  \cr
                C_{0} & \bullet &  \bullet &  \bullet  &  \bullet &  \bullet  } \hspace{1cm}\longrightarrow \hspace{1cm} \bordermatrix{ & & &  &  & \cr
                              C_{d+2}' & \bullet  &  *  &  *  &  *  &  *  \cr
                 C_{d+1}' & \bullet &   &\cr
                  C_{d}' & \bullet &  \bullet &  \bullet  &  \bullet &  \bullet \cr
                                & \vdots & \vdots & \ddots & \vdots & \vdots\cr
                                C_{1}' & \bullet &  \bullet &  \bullet  &  \bullet &  \bullet \cr
                                C_{0}' & \bullet &  \bullet &  \bullet  &  \bullet &  *  } $$

\begin{center} Fig 1  \end{center}

\noindent Now we extend the map $f:\t_N \to \R^d$ to a map $f':\t_{N'}\to\R^{d+1}$ by embedding $\R^d\subset \R^{d+1}$ where the last coordinate of $\R^{d+1}$ is zero:  Let $e_1,\cdots , e_{d+1}$ be the standard basis of $\R^{d+1}$ and $v_0,\cdots , v_N$ be the vertices of $\t_N$ and $v_0,\cdots , v_N, v_{N+1},\cdots , v_{N'}$ be the vertices of $\t_{N'}$. Let $x\in\t_{N'}$ and write $x=\sum\limits^{N'}_{0}\lambda_iv_i$ in terms of barycentric coordinates. Define $f' : \t_{N'}\to\R^{d+1}$ by
$$f'(x)=f(\lambda_0v_0+\cdots +\lambda_{N-1} v_{N-1} +(\lambda_N+\lambda_{N+1}+\cdots + \lambda_{N'})v_N) +(\lambda_N+\cdots +\lambda_{N'})e_{d+1}   $$
Then note that $f'$ naps additional vertices to $e_{d+1}$, is continuous and $f'|\t_N=f$. Since BMZ is true for the parameters $(d+1,\;r,\;\R^{d+1})$ and $C_0',C_1',\cdots , C_{d+2}'  $ is a special coloring of $\t_{N'}$, by assumption of reduction lemma there exists a family $\bar{\sigma_1}, \cdots , \bar{\sigma_r}$ of $r$ disjoint rainbow faces of $\t_{N'}$ such that
$$f'(\bar{\sigma_1})\cap\cdots\cap f'(\bar{\sigma_r})\neq \phi$$
Now let us consider the faces $\sigma_i=\bar{\sigma_i}\cap\t_N$ for $1\leq i\leq r$ of $\t_{N}$.  Note that at least one of the faces $\bar{\sigma_i}$ must be completely inside $\t_N$, otherwise we will have $r$ vertices of $\t_{N'}$ which are outside $\t_N$, but there are only  $r-1$ such vertices by definition, a contradiction.
Hence $\bar{\sigma_1}\cap \t_N=\bar{\sigma_1}=\sigma_1$ is a rainbow face of $\t_N$. Therefore $$f'(\bar{\sigma_1})\cap\cdots\cap f'(\bar{\sigma_r})\subset\R^d\subset\R^{d+1}$$
Now we also claim that each of $\sigma_1,..., \sigma_r$ is nonempty, because if some $\bar{\sigma_i},\;i\geq 2$, does not intersect $\t_N$, then $f'(\bar{\sigma_i})$ will lie in the positive part of  the last axis of $\R^{d+1}$, which means $f'(\bar{\sigma_1})\cap f'(\bar{\sigma_i})=\phi$, a contradiction. This implies that
$$ f(\sigma_i)\cap\cdots\cap f(\sigma_r)=f'(\bar{\sigma_1})\cap\cdots\cap f'(\bar{\sigma_r}) \neq\phi.$$
\end{proof}

\bn The following is a simple generalization of Proposition 2.1

{\Prop   Let us remove $q$ vertices $2\leq q\leq r-2$ from some color class, say $C_0$, and introduce one new color class $C_{d+2}$ of size $q$ so that we have the new coloring of $\t_N, $
\begin{equation}
 C_0,C_1, \cdots , C_{d+!}, C_{d+2},\;|C_0|=r-q+1,\;  \\|C_i|=i-1,\;|C_{d+1}|=1,\;|C_{d+2}|=q ~~\tag{2.1}\label{eq:2.1}
\end{equation}
Then the BMZ Theorem is true for the coloring  (2.1) of $\t_N$ also.}
\begin{proof}
We add $q$ vertices to the color class $C_0$ of the same color and add $r-q-1$ new vertices to the class $C_{d+2}$ of the same color and get the coloring
\begin{equation}
C_0',\; C_1', \cdots , C_{d+1}' ,~ C_{d+2}' ~~\tag{2.2}
\end{equation} such that $|C_i'|=r-1,~~\forall~i\neq d+1$ and $|C_{d+1}'|=1$
Now consider the simplex $\t_{N'},~N=(d+2)(r-1)$ containing $\t_N$ as the front face and define a map $f:\t_{N'}\to\R^{d+1}$ by mapping all additional vertices of $\t_{N'}$ onto $e_{d+1}$, and extending linearly. Then note that $\t_{N'}-\t_N$ has again only $r-1$ vertices.
Then using the fact that BMZ is true for  $(d+1, r,\R^{d+1})$ for the special coloring (2.1), we find by the same argument as in Prop 2.1 that BMZ is true for $f:\t_N\to\R^d$ for the coloring (2.1) also.
\end{proof}

\bn The following covers yet another possibility.

{\Prop  If we remove two vertices from two different color classes, say $C-0, C-1$, and  add a new color class $C_{d+2}$ having these two additional vertices, we will have the coloring
\begin{equation} C_0,~C_1,~C_2,\cdots , C_{d+1},~C_{d+2}~~\tag{3.1\label{eq:3.1}}\end{equation} such that $|C_0|=|C_1|=r-2,~|C_i|=r-1,~|C_{d+1}|=1,~|C_{d+2}|=2$ Then the BMZ is true for the coloring  (3.1)
}

\begin{proof}
Once again we embed the colorings (3.1) into the coloring \begin{equation} C_0',~C_1',~C_2',\cdots , C_{d+1}',~C_{d+2}'~~\tag{3.2\label{eq:3.2}}\end{equation} where $|C_i'|= r-1,~|C'_{d+1}|=1,~|C'_{d+2}|=r-1$ Then (3.2) is again a special coloring for  $(d+1, r,\R^{d+1})$. We define a map $f': \t_{N'}\to\R^{d+1},~N'=(d+2)(r-1)$ extending $f$ linearly as before. Then, again as in Proposition 2.1, this yields BMZ for $f:\t_N\to\R^d,~N=(d+1)(r-1)$.
\end{proof}

\bn {\bf Note :} It is clear from the above propositions that so long as  we remove any number of vertices from anywhere and  create only one new color class, the arguments of Proposition 2.1 will work and we will get BMZ for such general colorings. The argument will have to be changed when the number of new color classes so created exceeds one, as in the following.

{\Prop  Let us remove two vertices from some color class say $C_0$, and add two new color classes $C_{d+2},~C_{d+3}$ of size one each  such that we have  the coloring
\begin{equation} C_0,~C_1,~C_2,\cdots , C_{d+1},~C_{d+2},~C_{d+3}~~\tag{4.1\label{eq:4.1}}\end{equation}  with $|C_0|=r-3~|C_{d+1}|=|C_{d+2}|=|C_{d+3}|=1$  and $|C_i|=r-1$ for all other $i$'s. Then BMZ is true for the coloring (4.1) also.
}

\begin{proof}
Add two new vertices to $C_0$ of the same color as that of $C_0$, add $r-2$ new vertices to each of $C_{d+2},~C_{d+3}$ to get a new coloring of $\t_{N'},~N'=(d+3)(r-1)$
\begin{equation} C_0',~C_1',~C_2',\cdots , C_{d+1}',~C_{d+2}',~C_{d+3}'~~\tag{4.2\label{eq:4.2}}\end{equation} where $|C'_i|=r-1$ for each $i\neq d+1$
and $|C_{d+1}'|=1$. Then  (4.2) is a special coloring of $\t_{N'}$. Let $N=N_1=(d+1)(r-1),~N_2=(d+2)(r-1),~N'=N_3=(d+3)(r-1)$ and treat $\t_N$ as a front face of $\t_{N_2},~\t_{N_2} $ as a front face of $\t_{N_3}=\t_{N'}$. Let $e_i,~1\leq i\leq d+2$ be the standard basis of $\R^{d+2}$. Define a map $f':\t_{N'}\to\R^{d+2}$ as follows. Map one removed vertex of $C_0$ and last $r-2$ vertices of $C_{d+1}'$ to $e_{d+1}$, map other removed vertex of $C_0$ and $r-2$ vertices of $C_{d+2}'$ to $e_{d+2}$. Then as in Proposition  2.1, extend $f$ to $f': \t_{N'}\to\R^{d+2}$ linearly. We have $\t_{N_1}\subset\t_{N_2}\subset\t_{N_3},~\R^d\subset\R^{d+1}\subset\R^{d+2}$ and $\t_{N_3}-\t_{N_2},~\t_{N_2}-\t_{N_1}$ each have only $r-1$ vertices. Now using the BMZ for the special coloring (4.2) and the continuous map $f':\t_{N_3}\to\R^{d+2}$, we find $r$ disjoint rainbow faces $\sigma_1, \sigma_2,\cdots , \sigma_r$ of $\t_{N_3}$ such that
$$ f(\sigma_i)\cap\cdots\cap f(\sigma_r)\neq\phi.$$
Now arguing as in Proposition 2.1 again, we find that $\sigma_1\cap \Delta_{N_2},~\cdots , \sigma_r\cap \Delta_{N_2}$ are rainbow faces of $\t_{N_2}$ such that their $f'$-images have non empty intersection. Intersecting above rainbow faces with $\Delta_{N_1}$ and arguing once more, we find that $\sigma_1\cap \Delta_{N_1},~\cdots , \sigma_r\cap \Delta_{N_1}$ are  rainbow faces of $\t_{N_1}$ whose $f'$-images and therefore $f$-images have non empty intersection.
\end{proof}

\bn Now let us start with a general coloring \begin{equation} C_0,~C_1,~C_2,\cdots , C_{d+1},\cdots , C_{d+k}, \; k\geq 1~~\tag{5.1\label{eq:5.1}}\end{equation} of $\t_N,~N=(d+1)(r-1)$ where $|C_i|\leq r-1$. We introduce  a new color class $C_{d+k+1}=\{\phi\}$.

\bn We extend these color classes to the classes  \begin{equation} C_0',~C_1',~C_2',\cdots , C_{d+1}',~C_{d+2}',\cdots ,C_{d+k+1}'~~\tag{5.2\label{eq:5.2}}\end{equation}  such that $|C'_{d+k+1}|=1$
and $|C_i|=r-1,~~\forall~i\neq d+k+1$. Clearly $C_0',~C_1',~C_2',\cdots , C_{d+1}',~C_{d+2}',\cdots ,C_{d+k+1}'$ is a special coloring of $\t_{N'}, ~N'=(d+k+1)(r-1)$. We can now think that $C_{d+1}'$ in (5.2) has been obtained  by removing some vertices from $C_0\cup\cdots\cup C_d$ and giving them a new color. In this way we obtain a new  coloring $C_0',~C_1',~C_2',\cdots , C_{d+1}'$ such that $|C_0'\cup\cdots C_{d+2}'|=(d+2)(r-1)=N_2$. We can now assume that $N=N_1\subset N_2$. Similarly we get numbers $N_3,~N_4,\cdots ,N_{k+1}$ such that $\Delta_{N_1}\subset \Delta_{N_2}\subset\cdots \subset \Delta_{N_{k+1}}$ and $|\t_{N_{i+1}}-\t_{N_i}|=r-1$ for each $i$. Now we can also assume that $$\R^d\subset\R^{d+1}\subset\cdots\subset\R^{d+k+1}$$ where the last coordinate is taken as zero everywhere.  Let $e_1,\cdots , e_d,~e_{d+1},\cdots , e_{d+k+1}$ be the standard basis for $\R^{d+k+1}$. Then as in Proposition 2.4, we can define a map $f':\t_{N_{k+1}}\to\R^{d+k+1}$ such that $f'$ maps first $r-1$ additional vertices to $e_{d+1}$, second $r-1$ additional vertices to $e_{d+2}$ etc.
Then we extend $f$ to $f':\t_{N'}\to\R^{d+k+1}$ linearly. By the BMZ for parameters $(d+k+1,r,\R^{d+k+1})$ for the special coloring (5.2), we find that we have a family $\sigma_1, \sigma_2,\cdots , \sigma_r$ of disjoint rainbow faces of  $\t_{N'}$ such that
$$ f'(\sigma_i)\cap\cdots\cap f'(\sigma_r)\neq\phi.$$  Then arguing  as in Proposition 2.4 repeatedly, we find that $\sigma_1\cap \t_N,~\cdots , \sigma_r\cap \t_N$ is a family of $r$ disjoint rainbow faces of $\t_N$ whose $f$-images will have non empty intersection.

\bn The above analysis proves the following

{\Prop Let $d\geq 1,~r\geq 2$ prime and $N=(d+1)(r-1)$. Let $C_0,\cdots , C_{d+k},\;\; k\geq 1$ be a general coloring of $\t_N$ where $|C_i|\leq r-1$. Then assuming  BMZ for special colorings of $\t_{N'},~~N'=(d+k)(r-1)$,  we find that  BMZ is true for any general coloring  of $\t_{N}$.}

\section{Reduction Lemma for maps to  Manifolds}
The BMZ Theorem for a $d$-manifold $M$ is stated as follows \cite{bmz3}:
{\Thm {\bf (BMZ for Manifolds)}  Let $d\geq 1, r\geq 2$ a prime $N=(d+1)(r-1)$. Suppose $$C_0,~C_1,\cdots, C_m,~~m\geq d+1$$ is a general  coloring of the simplex $\t_N$ such that $|C_i|\leq r-1$. Then for any given continuous map $f: \t_N\to M$, where $M$ is a $d$-manifold, we can find a family of $r$ disjoint rainbow faces  $\sigma_1, \sigma_2,\cdots , \sigma_r$ of $\t_N$ such that $$\bigcap f(\sigma_i)\neq \phi .$$ }

\bn Let $\I=[0,1]$ denote the unit interval and $\bar{M}=M\times \I^k$

\bn {\bf Reduction Lemma :(\cite{bmz3} Lemma 2.1) } {\it Suppose the BMZ Theorem for manifolds is true for all special colorings for the parameters $(d+k, r, M\times \I^k),~k\geq 0$ . Then the BMZ Theorem for maps to manifolds is true for all general colorings  for the parameters $(d,r,M)$.}

\begin{proof}
The proof of above Lemma is exactly parallel to the proof of Proposition 2.5 proved earlier except that we replace $\R^{d+k}$ everywhere by $M\times \I^k$ and $\R^d \subset \R^{k+1}\subset ...\subset \R^{d+k}$ by $M\subset M\times \I .... \subset  M\times \I^k$.
\end{proof}

\bn The reduction lemma as stated by BMZ in ~\cite{bmz3} says the following

\bn {\bf Reduction Lemma : } {\it It suffices to consider a manifold $M$ which satisfies the property that  \begin{equation}(r-1)\dim M>r . \hbox{cohdim} ~M~~\tag{*} \end{equation}
where cohdim $M$ means the cohomological dimension of $M$.}

\bn Here the cohomological dimension means the largest integer $d$ such that the $d$-dimensional cohomology of $M$ does not vanish. If we use the classical definition of cohomological dimension which conforms to the geometric intuition, the inequality (*) may appear a bit surprising. In order to state it in conformity with classical definition of cohomological dimension, let us recall (see~\cite{nag} Appendix, Cohomological Dimension Theory,  and also ~\cite{deo}, ~\cite{deo1} for details).

{\Def Let $X$ be any topological space and $H^*(X,\Z)$ denote the $\check{C}$ech cohomology of $X$ with coefficients in $\Z$. Then the largest integer $n$ (if it exists) such that $H^q(X,A;\Z)=0~\forall~q>n$ and for all closed sets $A$ of $X$, is called the cohomological dimension of $X$. Otherwise we say that cohomological dimension of $X$ is infinite. We write it as cohdim$(X)$.}

\bn The role of the pair $(X,A)$ where $A$ is a closed set is very important. Since $\R$ is contractible, all the positive dimensional cohomologies of $\R^n$ vanish. This does not mean that that the cohdim $\R^n$ is zero. In fact, let us see that cohdim $(\R^n)=n.$  Let $A$ be any closed subset of $\R^n$. Then it is well known that $H^q(\R^n, A; \Z)= 0~\forall~q>n$ and for all closed sets $A$ of $\R^n$, which implies that cohdim$\R^n\leq n.$ On the other hand, it is easy to see $H^n(\R^n,\S^{n-1};\Z)\neq 0$ and so cohdim$\R^n\geq n.$ Thus we find that $\dim \R^n=\hbox{cohdim}~\R^n=n$. Also, for a paracompact Hausdorff space $X$, cohdim is a local property ~\cite{deo}, which means if $M$ is a paracompact manifold, then
$$\dim M=\hbox{cohdim} M. $$
Now going back to the BMZ Reduction Lemma for manifolds in \cite{bmz3}, observe that $M\times \I^k$ is a manifold of dimension $d+k ~(\dim M=d)$ and since $M\times \I^k$ is homotopically equivalent to $M,$ we find that $H^q(M\times\I^k)=0~\forall~q>d$. Thus the manifold $M''$ needed at the end of the proof of Lemma 2.1 in~\cite{bmz3} satisfies the condition that \begin{equation}(r-1)\dim M'' >r.d  ~~\tag{**} \end{equation} where $\dim M''=d+k$ and $H^q(M'';\Z)=0~\forall~q>d$. The manifold $M''= M\times\I^k$ clearly satisfies this property for $k$ suitably large.

The author is thankful to G$\ddot{u}$nter M. Ziegler for some valuable comments on this expository article.

\end{document}